\documentclass[11pt, reqno]{amsart}

\DeclareMathAccent{\mathring}{\mathalpha}{operators}{"17}

 \usepackage{color}
 
\newcommand{\mysection}[1]{\section{#1}
      \setcounter{equation}{0}}

\newtheorem{theorem}{Theorem}[section]
\newtheorem{lemma}[theorem]{Lemma}

\newtheorem{corollary}[theorem]{Corollary}

\theoremstyle{definition}
\newtheorem{assumption}{Assumption}[section]

\theoremstyle{remark}
\newtheorem{remark}{Remark}[section]

\newtheorem{example}{Example}[section]

\newcommand\cbrk{\text{$]$\kern-.15em$]$}}
\newcommand\opar{\text{\raise.2ex\hbox{${\scriptstyle | }$}\kern-.34em$($} }

 \makeatletter
 \def\dashint{%
 \operatorname%
 {\,\,\text{\bf--}\kern-.98em\DOTSI\intop\ilimits@\!\!}}
 \makeatother

\newcommand\frB{\mathfrak{B}}

\newcommand\bR{\mathbb{R}}

\newcommand\bZ{\mathbb{Z}}

\newcommand\cL{\mathcal{L}}

\newcommand\cS{\mathcal{S}}

\begin{document}

\title[Rate of convergence]
{Rate of convergence
 of difference approximations 
for uniformly nondegenerate   elliptic Bellman's
equations}

\author{N.V. Krylov}
\thanks{The  author was partially supported by
NSF Grant DMS-0653121}
\email{krylov@math.umn.edu}
\address{127 Vincent Hall, University of Minnesota,
 Minneapolis, MN, 55455}



\begin{abstract}
We show that the rate of convergence of 
solutions of 
finite-difference
  approximations for uniformly elliptic Bellman's
equations is of order at least $h^{2/3}$, where $h$ is the mesh size.
The equations are considered in smooth bounded domains.
\end{abstract}

\maketitle

The convergence of and error estimates for monotone and consistent
 approximations to
fully nonlinear, first-order PDEs were established a while ago by
 Crandall and Lions
\cite{CL} and Souganidis \cite{So}. 

The convergence 
of monotone and consistent
approximations for fully nonlinear, possibly degenerate 
second-order 
PDEs was first
proved in Barles and Souganidis \cite{BS}.  In a series 
of papers Kuo and 
Trudinger \cite{KT90,KT92,KT96}
  also looked in great detail at the issues of regularity and 
existence of such
approximations for uniformly elliptic equations.

There is also a probability part of the story, which
 started long before see Kushner \cite{Ku},
Kushner and Dupuis \cite{KD}, also see Pragarauskas \cite{Pr}.

However, in the above cited articles
apart from \cite{CL,So}, related to the first-order equations,
 no rate of convergence
was established. One can read more about the past development
of the subject in Barles and Jakobsen \cite{BJ07} and  the joint article
of Hongjie Dong and the author \cite{DK07}. We are going to discuss
only some results concerning{\em second-order\/} Bellman's equations,
which arise in many areas of mathematics such as  control
theory, differential geometry, and mathematical finance
(see Fleming and Soner \cite{FS06}, Krylov \cite{Kr77}) and which are most relevant
to the results of the present article.  

The first estimates of the rate of convergence
for second-order degenerate Bellman's equations appeared  in 1997
(see \cite{Kr97}). For equations with constant ``coefficients"
and arbitrary monotone  finite-difference approximations
it was proved in \cite{Kr97} that the rate of convergence
is $h^{1/3}$ if the error in approximating the true operators with
finite-difference ones is of order $h$ on {\em
three\/} times continuously differentiable functions. 
The order becomes
$h^{1/2}$ if the error in approximating the true operators with
finite-difference ones is of order $h^{2}$
on {\em four\/} times continuously differentiable functions
(see Remark 1.4 in \cite{Kr97}, which however contains
an arithmetical error albeit easily correctable. Also see
 Theorem 5.1 in \cite{Kr97}).
The main idea of \cite{Kr97} that the equation and its 
finite-difference approximation should play symmetric
roles is also used in the present article. The proofs
in \cite{Kr97} are purely analytical (in contrast with what
one can read in some papers mentioning \cite{Kr97}) 
even though sometimes probabilistic
{\em interpretation\/} of some statements are also given.
The next step was done in \cite{Kr99} where the 
so-called method of ``shaking  the coefficients''
was introduced to deal with the case of {\em degenerate\/}
parabolic Bellman's equations with variable coefficients.
The two sided error estimates were given: from the one side
of order $h^{1/21}$ and from the other $h^{1/3}$. Here $h$
(unlike in \cite{Kr97})
was naturally interpreted as the mesh size and
the approximating operators were assumed to approximate
the true operator with error of order $h$ on {\em three\/}
times continuously differentiable functions. Until now
it is not known whether or not it is possible to improve $1/21$
in the general setting of \cite{Kr99}.

However, what is possible is that one can get better
estimates if one uses some special approximations,
say providing the error 
  of order $h^{2}$ of approximating the main part of the
true operator
on {\em four\/} times continuously differentiable functions.
This was already mentioned in  Remark 1.4 of \cite{Kr97}
and  used by Barles and Jakobsen in \cite{BJ02}
to extend the results in \cite{Kr97} to  
equations with variable lower-order ``coefficients''.

One can also consider  special finite-difference
approximations, for instance, only containing
pure second-order differences in place of second-order
derivatives, when this $h^{2}$
approximating error is automatic. In such cases
the optimal rate $h^{1/2}$ was obtained in the joint
work of Hongjie Dong and the author \cite{DK07}
for
parabolic Bellman's equations with Lipschitz
coefficients in domains.
Both ideas of symmetry and   ``shaking the coefficients''
is used in \cite{DK07} as well as in \cite{DK05}.
In the paper by Hongjie Dong and the author \cite{DK05}
we consider among other things weakly nondegenerate
Bellman's equations with {\em constant\/} ``coefficients''
in the whole space
and obtain the rate of convergence $h$, where
$h$ is the mesh size. It may be tempting to say
that this result is an improvement of earlier
results, however it is just a better rate 
under different conditions.

 It is worth noting that
the set of equations satisfying the conditions
 in \cite{DK07} is smaller than the one in the papers by
Barles and Jakobsen \cite{BJ05,BJ07}, the results of which
obtained by using the theory of viscosity solutions
  guarantee the rate $h^{1/5}$. However, in the  examples
given  in
\cite{BJ05,BJ07}
of applications of the general scheme,  
for us to get the rate $h^{1/2}$,
we (only) need to add the requirement  that
the coefficients be twice differentiable
(see \cite{Kr08})
and in \cite{BJ05,BJ07} they are only assumed to be 
once differentiable and still the rate $h^{1/5}$
is guaranteed. One more point to be noted is that
in \cite{BJ07} parabolic equations are considered
with various types of approximation such as 
Crank-Nicholson and splitting-up schemes
 related to the time
derivative.

In the present paper we add two more restrictions
on the equations from \cite{DK07}: a) we require the coefficients
to be in $C^{1,1}(\bR^{d})$,
b) we require the equation to be uniformly nondegenerate.
In this case we obtain the rate of convergence
$h^{2/3}$, which was announced previously in
\cite{Kr07} for equations in the whole space.
 This time ``shaking'' is not needed as we explain in
 Remark \ref{remark 7.1.1}. It would be
very interesting to find a way to derive the results
of the present paper by using methods from
\cite{BJ02,BJ05,BJ07} or other methods
based on the theory of viscosity solutions.

 \mysection{Main results}
                                     \label{section 1.22.1}

 Let $A$ be a set and let
$$ 
a_{k}^{\alpha}=a^{\alpha}_{k}( x),\quad
 b_{k}^{\alpha}=b^{\alpha}_{k}( x),\quad
c^{\alpha}=c^{\alpha} ( x),\quad f^{\alpha}(x) 
$$
 be real-valued  
functions of $(\alpha, x)$ 
defined on $A \times\bR^{d}$ for $k=\pm1,...,\pm d_{1}$.
We assume that, for each $x=(x_{1},...,x_{d})
\in\bR^{d}$, these functions are bounded
 with respect
to $\alpha\in A$.
Also let some vectors 
$$
l_{k}=(l_{1k },...,l_{dk })\in B:=
\{x:|x|\leq1\}\subset\bR^{d}
$$ 
be defined
for $k=\pm1,...,\pm d_{1}$. This somewhat unusual
range of $k$ turns out to be convenient as we explain in
Remark \ref{remark 7.29.1}.
Consider the following Bellman's equation arising, for instance, in 
the theory
of controlled diffusion processes (see, for instance, 
\cite{FS06}, \cite{Kr77}):
\begin{equation}
                                              \label{1.21.1}
\sup_{\alpha\in A}\big[ 
 L^{\alpha}v(x) 
-c^{\alpha}( x)v(x)+f^{\alpha}( x) ]=0,
\end{equation}
where
\begin{equation}
                                                    \label{1.21.2}
L^{\alpha}v(x)=\sum_{i,j=1}^{d }a^{\alpha}_{, ij}(x)D_{ij}v(x)
+\sum_{i=1}^{d }b^{\alpha}_{ , i}(x)D_{i}v(x),
\end{equation}
$$
D_{i}=\frac{\partial}{\partial x_{j}},\quad D_{ij}=D_{i}D_{j},
$$
\begin{equation}
                                                    \label{1.21.3}
 a^{\alpha}_{ , ij}(x)=\sum_{|k|=1}^{d_{1}}
a^{\alpha}_{k}(x)l_{ik } l_{jk } ,\quad
b^{\alpha}_{ , i}(x)=
\sum_{|k|=1}^{d_{1}}b^{\alpha}_{k }(x)l_{ik } .
\end{equation}
As follows from the title we will be dealing with
uniformly elliptic operators and
it is well-known that for any uniformly elliptic operator of type
\eqref{1.21.2} there exist constant vectors $l_{k}$
(independent of $\alpha$) such that
representation \eqref{1.21.3} holds. 
In addition, the regularity properties of $a^{\alpha}_{ij}$
are inherited by $a^{\alpha}_{k}$ (see, for instance,
Theorem \ref{theorem 4.7.1} below).
This is also known for many
degenerate elliptic operators (see, for instance,
\cite{Kr08}). 

In what follows we adopt the summation convention
over all ``reasonable'' values of repeated indices.
Observe that
$$
a^{\alpha}_{, ij}D_{ij}v=a^{\alpha}_{k}D^{2}_{l_{k}}v,\quad
b^{\alpha}_{ , i}D_{i}v=b^{\alpha}_{k}D_{l_{k}}v,
$$
where $D^{2}_{l_{k}}v$ and $D_{l_{k}}v$ are the second
and the first derivatives of $v$ in the direction of
 $l_{k}$, that is 
$$
D^{2}_{l_{k}}v=l_{ik} l_{jk} D_{ij}v,\quad
D_{l_{k}}v=l_{ik} D_{i}v.
$$
Therefore equation \eqref{1.21.1} is rewritten as
\begin{equation}
                                              \label{7.22.1}
\sup_{\alpha\in A}\big[ 
 a^{\alpha}_{k}(x)D^{2}_{l_{k}}v(x)+
b^{\alpha}_{k}(x)D _{l_{k}}v(x)
-c^{\alpha}( x)v(x)+f^{\alpha}( x) ]=0,
\end{equation}
where, naturally, the summations with respect to $k$
are performed inside the sup sign.
We  approximate solutions of
\eqref{7.22.1} by solutions of finite-difference
equations obtained after replacing $D^{2}_{l_{k}}v$
and $D _{l_{k}}v$ with second- and first-order
differences, respectively, taken in the direction of $l_{k}$.

For any $x,\xi\in\bR^{d}$,
$h>0$, and function $\phi$ on $\bR^{d}$ introduce
$$
T_{h,\xi}\phi(x)=\phi(x+h\xi),\quad\delta_{h,\xi}
=\frac{ T_{h,\xi}-1}{h},\quad \Delta_{h,\xi}=\frac{T_{h,\xi}
-2+T_{h,-\xi}}{h^{2}}.
$$
When $\xi$ is one of the $l_{k}$'s we use the notation
$$
\delta_{h,k}=\delta_{h,l_{k}}, \quad
\Delta_{h,k}=\Delta_{h,l_{k}}, \quad
k=\pm1,...,\pm d_{1},
$$
in which the  finite difference approximation of \eqref{7.22.1}
is the following
\begin{equation}
                                                    \label{10.26.1}
\sup_{\alpha\in A} 
 [a_{k}^{\alpha}( x)\Delta_{h,k}v(x)
+b_{k}^{\alpha}( x)\delta_{h,k}v(x) 
-c^{\alpha}( x)v(x)+f^{\alpha}( x) ]=0.
\end{equation}

\begin{assumption}
                                         \label{assumption 12.17.2}
We are given a function $g$ on $\bR^{d}$
and two constants $\delta,
K\in(0,\infty)$ 
such that for all $\alpha\in A$ and $k$ on $\bR^{d}$ we have
$$
a^{\alpha}_{k}\geq\delta  ,\quad
c^{\alpha}\geq\delta
$$
and for $\phi=g,a^{\alpha}_{j},b^{\alpha}_{k},c^{\alpha}$,
$f^{\alpha}$,
$\alpha\in A$, and $ j,k\in\{\pm 1,...,\pm d_{1}\}$
we have that $\phi\in C^{1,1}(\bR^{d})$  and 
$$
\|\phi\|_{C^{1,1}(\bR^{d})}\leq K.
$$

\end{assumption}

\begin{assumption}
                                \label{assumption 12.17.1}
(i) We have $l_{k}=-l_{-k}$, 
$a^{\alpha}_{k}=a^{\alpha}_{-k}$, $k=\pm1,...,\pm d_{1}$.

(ii) There exists an integer $1\leq d_{0}\leq d_{1}$ such that
for   
$$
\Lambda :=\{l_{k},k=\pm1,...,\pm d_{1}\} ,\quad
  \cL :=\{ l_{\pm1},..., l_{\pm d_{0}} \} 
$$
we have $0\in\cL$ and
$$
\cL+\cL\supset \Lambda\supset\{l'+l'':l',l''\in\cL,l'\ne l'' 
 \}.
$$

(iii) The coordinates of $l_{k}$ are rational numbers
and $\text{Span}\,\Lambda=\bR^{d}$.
\end{assumption}
 
To justify Assumptions \ref{assumption 12.17.2}
 and \ref{assumption 12.17.1} we remind the reader
Theorem 3.1 of \cite{Kr11} in which $\bZ=\{0,\pm 1,\pm 2,...\}$
and
$\cS_{\delta_{1}}$, $\delta_{1}>0$, is the set  
of symmetric $d\times d$-matrices $a$
 such that for any $\xi\in\bR^{d}$
$$
\delta_{1}|\xi|^{2}\leq \langle a\xi,\xi\rangle
\leq\delta_{1}^{-1}|\xi|^{2} 
$$
($\langle \cdot,\cdot\rangle$ stands for the scalar
product).
\begin{theorem}
                                  \label{theorem 4.7.1}
There exists a set $\{l_{1},...,l_{n}\}
\subset \bZ^{d}$     such that
for any its extension $\{l_{1},...,l_{m}\}$, $m\geq n$,
  there exist
real-analytic functions $\lambda_{1}(a),...,
\lambda_{m}(a)$ on $\cS_{\delta_{1}}$ possessing there
the following properties:
$$
a\equiv\sum_{k=1}^{m}\lambda_{k}(a)l_{k}l_{k}^{*},
\quad \lambda_{k}(a)\geq\delta 
,\quad \forall k,
$$
where the constant $\delta >0$.

\end{theorem}

\begin{remark}
                                \label{remark 7.29.1}
Theorem \ref{theorem 4.7.1} implies that
any uniformly nondegenerate equation
of type \eqref{1.21.1} can be written as 
\eqref{7.22.1} with the coefficients
satisfying Assumption \ref{assumption 12.17.2}
as long as the coefficients in \eqref{1.21.1}
are bounded and twice continuously differentiable
with $C^{2}$-norm controlled by a constant
independent of $\alpha$ and $c$ is uniformly 
bounded away from zero. If we take
$ \{l_{1},...,l_{n}\}$ from Theorem \ref{theorem 4.7.1}
and define $\cL=\{0,\pm l_{1},...,\pm l_{n}\}$ and
$$
\Lambda=\{l'+l'':l',l''\in\cL,l'\ne l'' 
 \},
$$
then Assumption \ref{assumption 12.17.1}(ii)
will be obviously satisfied and owing to Theorem
\ref{theorem 4.7.1} Assumption \ref{assumption 12.17.2}
will still be preserved. Assumption \ref{assumption 12.17.1}(i)
is of a technical nature and easily satisfied just by redefining
$a^{\alpha}_{k}$ if necessary
 which is possible since the above $\Lambda$
is symmetric and  $D^{2}_{l }=
D^{2}_{-l }$. One could exclude 
Assumption \ref{assumption 12.17.1}(i) on the expense
of more complicated formulation
of Assumption \ref{assumption 12.17.1}(ii), which, by the way,
is needed in order to apply the results
of \cite{Kr12.2} about interior estimates
of second-order differences of approximate solutions.
It is also worth saying that by assumption our
$l_{k}$ have lengths $\leq1$. Therefore one should
normalize those $l_{k}$'s from Theorem \ref{theorem 4.7.1},
which are different from zero,
absorbing the extra factors in $\lambda_{k}(a)$.

Also observe that Assumption \ref{assumption 12.17.1} 
requires that $l_{k}=0$ for some $k$.
 Including the zero vector in $\Lambda$
turns out to be  convenient from technical point of view. The coefficient
$a^{\alpha}_{k}$ corresponding to this vector can be set to equal,
say 1, because the corresponding finite-difference operator is just zero.

Finally, talking about our assumptions we point out that
  Assumptions  \ref{assumption 12.17.1}(ii)(iii) 
are used in Section \ref{section 7.4.1}
when we apply some results of
\cite{Kr12.1} and \cite{Kr12.2} to finite difference
equations in domains.
\end{remark}

We suppose that we are given a $\psi\in C^{2}(\bR^{d})$  
such that
$$
\Omega=\{x:\psi(x)>0\}
$$
is a bounded domain and $| D \psi|\geq1$ on $\partial \Omega$.
Equation \eqref{7.22.1} is considered in $\Omega$
with $v$ subject to the   boundary condition $v=g$
on $\partial\Omega$.

Introduce
$$
\Omega_{h}=\{x\in \Omega:x+hB\subset \Omega\},\quad
\partial_{h}\Omega= \bR^{d}\setminus \Omega_{h}.
$$ 
 Observe that $\partial_{h}\Omega$ contains points
which are very far from $\Omega$. This turns out to be convenient
in our constructions.

Here are our main results in which the above assumptions
are supposed to hold and
$$
\rho(x)=\text{dist}\,(x,\Omega^{c}),\quad\text{where}\quad
\Omega^{c}=\bR^{d}\setminus\Omega.
$$
Naturally, $\rho=0$ on $\Omega^{c}$.
\begin{theorem}
                                        \label{theorem 7.22.1}
There are  constants $N\in[0,\infty)$ and $h_{0}>0$ such that
for all  $h\in(0,h_{0}]$
  the following is true:

(i) Equation \eqref{10.26.1} in $\Omega_{h}$ with  
boundary condition $v=g$ on $\partial_{h}\Omega$ has a unique
bounded Borel solution $v _{h}$.

(ii) On $\bR^{d}$
\begin{equation}
                                                    \label{6.25.3}
\rho^{-1}|v _{h}-g |,\quad|\delta_{h,i}v _{h}|,\quad
  (\rho-6h)|\delta_{h,i}\delta_{h,j} v_{h}|\leq N
\end{equation}
for any $i,j$ and for any $x,y\in\bR^{d}$
\begin{equation}
                                                    \label{7.29.1}
|v _{h}(x)-v _{h}(y)|\leq N(|x-y|+h).
\end{equation}

\end{theorem}

\begin{theorem}
                                        \label{theorem 7.22.2}

 There exists a unique $v \in C^{1,1}_{loc}(\Omega)\cap
 C^{0,1}(\bR^{d})$
satisfying 
equation \eqref{7.22.1} in $\Omega$ (a.e.) and
equal $g$ on $\Omega^{c}$. Furthermore,  
 $\rho|D^{2}v |\leq N$ in $\Omega$ (a.e.),
where $N$ is a constant.

\end{theorem}
This theorem is proved in exactly the same way in which 
Theorem 8.7 of \cite{Kr11} is proved. On this way
one uses Theorem \ref{theorem 7.22.1}, 
 the fact that the derivatives of $v$ are weak limits
of finite differences of $v_{h}$ as $h\downarrow0$
(see the proof of Theorem 8.7 of \cite{Kr11}), and the fact 
that 
  there are sufficiently many
  second order derivatives in directions
of $l_{i}$, $l_{j}$ to conclude from their boundedness
that the Hessian of $v$ is bounded.

\begin{remark}
                                          \label{remark 6.26.1}
In the theory of fully nonlinear elliptic equations
much more general results than Theorem \ref{theorem 7.22.2}
 under much weaker
conditions are known. For instance, it is proved in \cite{Sa}
that if $\beta\in(0,1)$ is sufficiently small
and $a^{\alpha},b^{\alpha},c^{\alpha},f^{\alpha}\in C^{\beta}(\bar{\Omega})$
with $C^{\beta}(\bar{\Omega})$-norms bounded by a constant
independent of $\alpha$, then $v\in C^{2+\beta}_{loc}(Q)\cap C^{0,1}(\bar{Q})$.

The results in \cite{Sa} and other classical texts on
the theory of fully nonlinear elliptic equations are obtained
on the basis of very deep facts, using very sophisticated 
and beautiful techniques,
and require a series of long arguments in the end of which
the reader learns a lot of various facts from the theory
of PDEs and functional analysis.
In contrast, our Theorem \ref{theorem 7.22.2} is obtained
on the sole basis of the discrete maximum principle
combined with elementary albeit  quite long computations
(see \cite{Kr12.1} and \cite{Kr12.2}) involving
discrete versions of Bernstein's method.

\end{remark}

\begin{theorem}
                                        \label{theorem 7.22.3}
There are  constants $N\in[0,\infty),h_{0}>0$ such that
for all  $h\in(0,h_{0}]$ we have $|v_{h}-v|\leq Nh^{2/3}$
on $\bR^{d}$. 
\end{theorem}

\begin{example}
Consider the following uniformly nondegenerate analog
of the Monge-Amp\`ere equation
\begin{equation}
                                                  \label{7.5.1}
{\rm det}\,(-D_{ij}v-\gamma^{2}\delta^{ij}\Delta v)=
(f_{+})^{d},
\end{equation}
where $\gamma>0$ is a parameter and $f\in C^{1,1}(\bR^{d})$.
It is well known (see, for instance, \cite{Kr72}) 
that equation \eqref{7.5.1} supplied with the requirement
that the matrix $( D_{ij}v+\gamma^{2}\delta^{ij}\Delta v)$
be negative definite is equivalent to the following
single Bellman's equation:
\begin{equation}
                                                  \label{7.5.2}
\sup_{\substack{a=a^{*}\geq0,\\{\rm trace}\, a=1}}
[(a_{ij}+\gamma^{2}\delta^{ij})D_{ij}v+
({\rm det}^{1/d}a)fd]=0.
\end{equation}
The above theory is applicable to \eqref{7.5.2} and
shows that we can approximate the solutions of
\eqref{7.5.1} with a $C^{1,1}$ boundary
condition satisfying $( D_{ij}v+\gamma^{2}\delta^{ij}\Delta v)
\leq0$ with solutions of corresponding finite-difference
equations. It is to be noted however that
the number $d_{1}$ related to the number
of directions needed to write \eqref{7.5.2}
in the form \eqref{7.22.1} will depend on $\gamma$ and will go
to infinity as $\gamma\to0$.

Also observe that one can prove that if $\Omega$
is a strictly convex domain, and $v^{\gamma}$
are solutions of \eqref{7.5.2} with zero boundary condition,
 then $|v^{\gamma}-v^{0}|
\leq N\gamma$, where $v^{0}$ is a (generalized or viscosity)
 solution of
\eqref{7.5.1} with $\gamma=0$.
\end{example}

The rest of the article is organized as follows.
In Section \ref{section 7.25.1} we prove
\eqref{7.29.1} assuming \eqref{6.25.3} for equations
more general than Bellman's equations.
Section \ref{section 7.4.1} contains two
results from \cite{Kr12.1} and 
 \cite{Kr12.2} needed to prove 
\eqref{6.25.3} in Section \ref{section 7.5.1}.
In the final rather long Section \ref{section 7.5.2}
we present the proof of Theorem \ref{theorem 7.22.3}.

We use $N$ to denote various constants which may change from one
appearance to another. Sometimes we specify what they
are depending on or independent of. However sometimes 
we do not do that.
In these situations it is understood that they
are independent of anything which is allowed to change
like $x,h,\varepsilon$....

\mysection{On the Lipschitz continuity
of $v_{h}$}

                                            \label{section 7.25.1}

The setting and notation in this section
are different from the ones in Section \ref{section 1.22.1}. 
Let
$$
\Lambda:=\{ \ell_{k};k=  1,...,  d_{1}\} 
$$
be a symmetric subset of   $B=\{x:|x|<1\}$.   

 Let $H(p,x,u)$ be a real-valued function
given for

$$
 p\in\bR,\quad x\in
\bR^{d},\quad u=(u',u''),\quad u'=(u'_{0},u'_{1},
...,u'_{d_{1}}),
\quad u''=
(u''_{1},...,u''_{d_{1}}).
$$

Fix two constants $\delta\in(0,1],K\in[0,\infty)$.

\begin{assumption}
                                         \label{assumption 6.11.1}
The function $H(p,x,u)$ is locally Lipschitz continuous with respect 
to $(p,u)$. Furthermore, at all points of differentiability 
of $H$ with respect to $(p,u)$
we have
$$
\delta \leq H_{u''_{j}} \leq K ,
\quad j=1,...,d_{1},\quad H_{u'_{0}}\leq-\delta,
$$
$$
|H_{u'_{j}}|\leq K,\quad j=0,...,d_{1},\quad
|H_{p}(p,x,u)|\leq K(1+|u|) .
$$

\end{assumption}

Introduce
$$
H(x,u)=H(0,x,u).
$$

The following lemma is an easy consequence of the mean value
theorem in the form of Hadamard.
\begin{lemma}
                                       \label{lemma 6.17.1}
For any $p\in\bR $, $x\in\bR^{d}$, and $u,v\in\bR^{2d_{1}+1}$
there exist  numbers $a_{1}$,..., $a_{d_{1}}$,
$b_{1}$,..., $b_{d_{1}}$,
 $c$, and $f$ such that
$$
H(p,x,u)-H(x,v)=H(p,x,u)-H(p,x,v)+f
$$
$$
=a_{i}(u''_{i}-v''_{i})+b_{i}
(u' _{i}-v' _{i})+c(u'_{0}-v'_{0})+f,
$$
\begin{equation}
                                                  \label{6.17.1}
\delta\leq a_{i}\leq K,\quad -K\leq c\leq-\delta,\quad
|b_{i}|\leq K,\quad |f|\leq K|p|(1+ |v|)
\end{equation}
for all $i\in\{1,...,d_{1}\} $.
\end{lemma}

For any function $v$ and $h>0$ define
$$
H_{h}[v](x)=H( x,v (x),\delta_{h }v (x),
\Delta_{h}v (x)),
$$
where
$$
\delta_{h }v =(\delta_{h,1}v ,...,
\delta_{h,d_{1} }v ),
\quad \Delta_{h}v =(\Delta_{h,1} v ,...,
\Delta_{h,d_{1}} v ),
$$
$$
\delta_{h,k}=\delta_{h,\ell_{k}},\quad
\Delta_{h,k}=\Delta_{h,\ell_{k}}.
$$
 
\begin{theorem}
                                            \label{theorem 7.12.1}
Let   $D$ be a domain in $\bR^{d}$ and let  $h\in(0,\delta/K]$. 
Then there exists a constant $N_{0}$ depending only on
$\delta,K$, and $d_{1}$ such that for any two functions
$v'$ and $v''$ given on $\bR^{d}$ we have that in $D_{h}$
\begin{equation}
                                                       \label{5.15.6}
N_{0}(v' -v'' )+
h^{2}\big(H_{h}[v'] -H_{h}[v''] \big)
\leq
(N_{0}- \delta  h^{2} )
\sup_{D}(v'-v'')_{+},
\end{equation}
\begin{equation}
                                                       \label{7.12.2}
\sup_{\bR^{d}}(v'-v'')_{+}\leq\delta^{-1}
\sup_{D_{h}}(H_{h}[v''] -H_{h}[v'] \big)_{+}
+\sup_{\partial_{h} D}(v'-v'')_{+}.
\end{equation}

Furthermore, if $H( x,0)$ is bounded and
we are given a bounded function $g$ on $\bR^{d}$, then the
  equation 
\begin{equation}
                                              \label{2.25.3}
H_{h}[v]=0  
\end{equation} 
in $ D_{h} $ 
with boundary condition
$v=g$ on $\partial_{h}D$ (in case $\partial_{h}D\ne\emptyset$)
has a unique bounded solution.
This solution is Borel measurable in $x$ if $g$
 and $H( x,u)$ are Borel for each $ u$.
\end{theorem}

Proof. As is usual for monotone finite-difference equations
in the proof of solvability of \eqref{2.25.3}
we rely on 
the method called Jacobi iteration in \cite{KT92} and
the Banach fixed point theorem
in the space $\frB$ of bounded functions on
$\bR^{d}$ provided with the sup norm. 
We will follow an argument from \cite{Kr11}
where the situation is somewhat different.

First we deal with \eqref{5.15.6}.
By Lemma \ref{lemma 6.17.1},
 for   $w=v'-v''$
we have  on $\bR^{d}$ that
$$
h^{2}\big(H_{h}[v'](x)-H_{h}[v''](x)\big)
=(ch^{2}-h \sum_{j=1}^{d_{1}}b_{j}-2 \sum_{j=1}^{d_{1}}
a_{j})w(x) 
$$
$$
+h b_{j} w (x+h\ell_{j}) + 
a_{j} [w(x+h\ell_{j})  +
v (x-h\ell_{j})].
$$

Since $a_{j},b_{j},-c\leq K $, and $h\in(0,\delta/K]$, 
there is a constant $N_{0}$  depending only on
$\delta,K$, and $d_{1}$ such that
$$
N_{0}+ch^{2}-h \sum_{j=1}^{d_{1}}b_{j}-2 \sum_{j=1}^{d_{1}}a_{j}\geq0.
$$
Furthermore, 
$$
h b_{j} w(x+h\ell_{j}) + 
a_{j} [w(x+h\ell_{j})  +
w(x-h\ell_{j})]= (hb_{j}+a_{j})w (x+h\ell_{j})
+  a_{j} w(x-h\ell_{j})
$$
and all the coefficients on the right are nonnegative
in light of the fact that
 $h\in(0,\delta/K]$ and $a_{j}\geq\delta$. It follows that
in $D_{h}$
the left-hand side of \eqref{5.15.6} is less than
$$
 \big(N_{0}+ch^{2}-h \sum_{j=1}^{d_{1}}b_{j}-2 \sum_{j=1}^{d_{1}}a_{j}
\big)\sup_{D} w_{+}
$$
$$
+\sum_{j=1}^{d_{1}}(hb_{j}+a_{j}) \sup_{D}w_{+}
 + \sum_{j=1}^{d_{1}} a_{j}  \sup_{D}w_{+}
$$
$$
 =(N_{0}+ch^{2})
\sup_{D}(v'-v'')_{+}\leq
(N_{0}- \delta  h^{2} )
\sup_{D}(v'-v'')_{+},
$$
which proves \eqref{5.15.6}. 

While proving \eqref{7.12.2} we may assume that
$$
\sup_{\bR^{d}}(v'-v'')_{+}>\sup_{\partial_{h}D}(v'-v'')_{+}.
$$
In that case 
$$
\sup_{\bR^{d}}(v'-v'')_{+}=\sup_{D_{h}}(v'-v'')_{+}
$$
and \eqref{5.15.6} implies
$$
N_{0}(v' -v'' )_{+}
\leq
(N_{0}- \delta  h^{2} )
\sup_{D_{h}}(v'-v'')_{+}+
h^{2}\big(H_{h}[v''] -H_{h}[v'] \big)_{+}.
$$
After that
it only remains to take the sups of both sides over $D_{h}$
and collect like terms.
 
Now we prove the second assertion of  the theorem.
Observe that, as is easy to see, equation
\eqref{2.25.3} in $D_{h}$
with the boundary condition $v=g$ on $\partial_{h}D$
 is equivalent to the following single
equation:
\begin{equation}
                                                       \label{5.15.1}
v(x)=\big[N_{0}^{-1}h^{2}
H_{h}[v]( x)
+v(x)\big]I_{D_{h}}(x)+g(x)I_{\partial_{h}D}(x).
\end{equation}

Then   introduce an operator 
$T_{h}:v\to T_{h}v$
  by
$$
T_{h}v(x):=\big[N_{0}^{-1}h^{2}H_{h}[v]( x)
+v(x)\big]I_{D_{h}}(x)+g(x)I_{\partial_{h}D}(x).
$$
Since $H[0]$ is bounded by assumption, $T0$
is bounded.
By \eqref{5.15.6}  we have that
$T_{h}v=T_{h}v-T_{h} 0  +T_{h} 0 $ is bounded if
$v$ is bounded, so that $T_{h}$ is an operator in $\frB$,
and moreover $T_{h} $
is a contraction operator 
in $\frB$   with contraction
constant   less than $1-N_{0}^{-1} \delta h^{2} 
<1$. By the Banach fixed point theorem, 
equation \eqref{5.15.1} has a unique bounded solution.
This solution can be obtained as the limit of $(T_{h})^{n}0$
as $n\to\infty$. Furthermore, $T_{h}$ maps Borel functions
into Borel measurable ones if $H( x,u)$ is Borel with respect to $x$.
Therefore, under this condition, given that $g$ is Borel,
all functions 
$ T_{h} ^{n}0$ are Borel measurable, so that the solution
$v$ is also Borel measurable.
The theorem is proved.

By taking first $v'=v$, $v''=0$ and then $v'=0$, $v''=v$
we obtain from \eqref{7.12.2} the following. 

\begin{corollary}
                                       \label{corollary 7.12.1}
If $v$ is a solution of \eqref{2.25.3} with  
boundary condition $v=g$ on $\partial_{h}D$, then
$$
\sup_{\bR^{d}}v_{+}\leq\delta^{-1}\sup_{D_{h}}
\big(H_{h}[0])_{+}+\sup_{\partial_{h}D}g_{+},\quad
\sup_{\bR^{d}}v_{-}\leq\delta^{-1}\sup_{D_{h}}
\big(H_{h}[0])_{-}+\sup_{\partial_{h}D}g_{-}.
$$
\end{corollary}
 
We improve these estimates
in Lemma \ref{lemma 7.4.1}
 for more restricted range of~$h$.
We need a version of   Lemma 8.5 of \cite{Kr11}.
Let $\Omega$ be the set introduced in 
Section~\ref{section 1.22.1}.
\begin{lemma}
                                        \label{lemma 7.7.1}
Assume that ${\rm Span}\,(\Lambda)=\bR^{d}$.
Then there exist $h_{0}\in(0,\infty)$
and a nonnegative function $\Psi\in C^{2}(\bar{\Omega})$
 such that   $\Psi/\rho$ and $\rho/\Psi$ are bounded in $\Omega$,
for   $h\in(0,h_{0}]$ on $\Omega_{h}$ we have
\begin{equation}
                                           \label{7.3.1}
a_{j}\Delta_{h,j}\Psi + K
\sum_{j=1}^{m}|\delta_{h,j}\Psi |\leq-1,
\end{equation}
and such that there exist   constants $M,\mu\geq1$ for which
the function $\Phi:=M\Psi\ln(M/\Psi)$
satisfies
\begin{equation}
                                           \label{5.15.4}
a_{j}\Delta_{h,j}\Phi  
+ K\sum_{j=1}^{m}|\delta_{h,j}\Phi  |
\leq -\rho^{-1} 
\end{equation}
in $\Omega_{\mu h}$, whenever $\delta\leq a_{j}\leq K$. 
\end{lemma}

Below by $h_{0}$ we mean a constant in $(0,\delta/K]$
for which the statement of Lemma 
\ref{lemma 7.7.1} is true.

\begin{lemma}  
                                       \label{lemma 7.4.1}
         Assume that ${\rm Span}\,(\Lambda)=\bR^{d}$.
Let $D\subset \Omega$, $h\in(0,h_{0}]$ and let $v$ satisfy
\eqref{2.25.3} in $D_{h}$ and $v=g$ in $\partial_{h}D$,
where $g\in C^{1,1}(\bR^{d})$.
Then $|v-g |\leq N \rho $ on $\bR^{d}$, where the constant
$N$ is independent of $h$.
\end{lemma}

Proof. It suffices to prove
that $|v -g|\leq N \Psi $ on $\bR^{d}$ if we continue $\Psi$ outside
$\Omega$ as zero.
 In order to do this
we observe that by Lemma \ref{lemma 6.17.1} with $p=0$
at each point of $D_{h}$
$$
H_{h}[g+N\Psi]=H_{h}[g+N\Psi]-H_{h}[0]+H_{h}[0]
$$
$$
=N(a_{j}\Delta_{h,i}\Psi +  b_{j}\delta_{h,j}\Psi 
+ c\Psi)+ 
a_{j}\Delta_{h,i}g +  b_{j}\delta_{h,j}g 
+ cg+H(\cdot,0),
$$
where $a_{j},b_{j},c $ are some numbers satisfying
\eqref{6.17.1}. Owing to \eqref{7.3.1} the last 
expression is negative if we take $N$ large enough.
This implies that $v -g\leq N \Psi $
 by the maximum principle (that is  by \eqref{7.12.2})
and by the fact that $v =g$ and $\Psi\geq0$
 in $\partial_{h}D$. The estimate $v -g\geq -N\Psi $
is proved similarly. The lemma is proved.

In the rest of the section
we 
assume that 
$$
{\rm Span}\,(\Lambda)=\bR^{d},
$$
fix a  constant  $\mu\geq1$, 
and take a constant  $ \nu\in[0,\infty)$. Define
$$
\kappa=\mu\vee(2\nu)
$$
and let $D$ be an open subset of $\Omega_{\kappa h}$. For each $p\in\bR$
and $h>0$
introduce an operator $H^{p}_{h}$ acting on functions $v$ given
on $\bR^{d}$ by the formula
$$
H^{p}_{h}[v](x)=H(p,x,v (x),\delta_{h }v (x),
\Delta_{h}v (x)).
$$

\begin{theorem}
                                   \label{theorem 6.5.1}
Assume that for each $h\in(0,h_{0}]$ and
each $p$ we are given a bounded function  $v^{p}_{h}$ on
$\bR^{d}$ which satisfies the equation
\begin{equation}
                                             \label{6.11.1}
H^{p}_{h}[v^{p}_{h}]=0
\end{equation}
in $D_{h}$. Introduce
$$
  v_{h}=v^{0}_{h} 
$$
and assume that there is a constant $N_{0}\in[0,\infty)$
  such that
for     $h\in(0,h_{0}]$ and $x\in D_{h}$ 
 we have
$$
|v _{h}(x)|,\quad
|\delta_{h,i} v _{h}(x)| ,\quad (\rho (x)-\nu h)
|\Delta_{h,i} v _{h}(x)|\leq N_{0}
$$
for all $i$. Then for  $h\in(0,h_{0}]$ 
\begin{equation}
                                                \label{3.26.06}
|v^{p}_{h}(x)-v_{h}(x)|\leq N (|p|+  \sup_{ \partial_{ h}D}
|v^{p}_{h}-v_{h}| )
\end{equation}
for any $x \in\bR^{d}$ and $p\in\bR$,
 where the constant $N$ is independent
of $x,p$, and~$h$.

\end{theorem}

Proof. 
Take a constant $N_{1}$ to be specified later and
introduce 
$$
w_{h}=v_{h}+N_{1}p\Phi.
$$
 By Lemma \ref{lemma 6.17.1}
for each $x$
$$
H^{p}_{h}[ w_{h}](x) 
 -H[v_{h}] (x) 
$$
\begin{equation}
                                           \label{6.20.1}
=N_{1}pa_{j}\Delta_{h,i}\Phi(x)+N_{1}pb_{j}\delta_{h,j}\Phi(x)
+N_{1}pc\Phi(x)+f,
\end{equation}
where $a_{j},b_{j},c $ are some numbers satisfying
\eqref{6.17.1}  and $f$ is such that
$$
|f|\leq K|p| (1+ \sum_{j}(|\Delta_{h,j}v_{h}(x)|
+|\delta_{h,j}v_{h}(x)|)+ |v_{h}(x)|).
$$
 Observe that the second term on the left in \eqref{6.20.1}
vanishes in $D_{h}$. Furthermore, $0\leq(\rho-\nu h)^{-1}
\leq 2\rho^{-1}$ if $\rho\geq2\nu h$, so that, owing to the fact that
$\kappa\geq2\nu$ and $D_{ h}\subset\Omega_{\kappa h}$, in $D_{ h}$
we obtain
$$
|f|\leq N_{2}|p| \rho^{-1},
$$
where $N_{2}$ is independent of $h$ and $p$.
 Finally, we set $N_{1}=N_{2}$, take into account
\eqref{5.15.4} and the fact that $D_{h}\subset
\Omega_{\mu h}$ (since $\kappa\geq\mu$), and we 
conclude from \eqref{6.20.1} that
\begin{equation}
                                                \label{7.21.2}
H^{p}_{h}[ w_{h}]\leq0
\end{equation}
 in $D_{ h}$ if $p\geq0$. Upon comparing this with
 \eqref{6.11.1}
and using \eqref{7.12.2} we get that for any $x$
$$
 v^{p}_{h}(x )\leq w_{h}(x)
+\sup_{\partial_{ h}D} (v^{p}_{h}-w_{h})_{+}
$$
\begin{equation}
                                                \label{7.21.3}
\leq v_{h}(x)+N p +\sup_{\partial_{ h}D} (v^{p}_{h}-v_{h})_{+},
\end{equation}
where $N$ is independent of $x,h$ and $p$. 

If $p\leq0$ the inequality in \eqref{7.21.2} is
reversed and one gets
$$
v_{h}-Np\leq
w_{h}\leq v^{p}_{h}+\sup_{\partial_{ h}D} (w_{h}-v^{p}_{h})_{+}
\leq v^{p}_{h}+\sup_{\partial_{ h}D} (v_{h}-v^{p}_{h})_{+}.
$$
By combining this with \eqref{7.21.3} we come to
\eqref{3.26.06} and the theorem is proved.

\begin{corollary}
                                       \label{corollary 7.21.1}
Assume that $H(x,u)$ is locally Lipschitz continuous with
respect to $(x,u)$ and at all points of its differentiability
with respect to $(x,u)$
$$
|H_{x_{i}}(x,u)|\leq K(1+|u|)
$$
for all $i$. Let $H(x,0)$ be bounded on $\bR^{d}$.
For $h\in(0,\delta/K]$ denote by $v_{h}$ a unique
bounded solution of \eqref{2.25.3} in $\Omega_{h}$
with   boundary condition  $v=g$ on $\partial_{h}\Omega$
(see Theorem \ref{theorem 7.12.1}), where $g
\in C^{0,1}(\bR^{d})$.
 
Finally, assume that
there is a constant  $N_{0}\in[0,\infty)$
such that
for   $h\in(0,h_{0}]$ in $\Omega_{h}$ we have
$$
|\delta_{h,i} v_{h}(x)| ,\quad (\rho (x)-\nu h)
|\Delta_{h,i} v_{h}(x)|\leq N_{0}
$$
for all $i$.

 Then for all $h\in(0,h_{0}]$ and  $x,y\in \bR^{d}$
\begin{equation}
                                                \label{3.26.6}
|v_{h}(x)-v_{h}(y)|\leq N (|x-y|+ h ),
\end{equation}
  where the constant  $N$ is independent
of $x,y$, and $h$.

\end{corollary}

Proof. 
If $|x-y|\geq h$ one can split the straight
segment between
$x$ and $y$ into adjacent pieces of length $h$ combined with
a remaining one of length less than $h$. This shows that
we need only  prove \eqref{3.26.6} for $|x-y|\leq h$.

Then fix a unit vector $l\in\bR^{d}$ and for $p\in\bR $
redefine $H$ if necessary by setting
$$
H(p,x,u)=H(x+p l,u)\quad\text{if}\quad |p|\leq h,
$$
$$
H(p,x,u)=H(x+ l \,{\rm sign}\, p,u)\quad\text{if}\quad |p|\geq h,
$$
that is
$$
H(p,x,u)=H(x+  l \phi(p),u),
$$
where $\phi(p)=(-h)\vee p\wedge h$. Observe that the
function
$$
v^{p}_{h}(x):=v_{h}(x+ l \phi(p))
$$
satisfies \eqref{6.11.1} (with new $H$) in $\Omega_{2h}
\supset\Omega_{(\kappa+1) h}$, where the inclusion
follows from the fact that $\kappa\geq\mu\geq1$. 
By Lemma \ref{lemma 7.4.1} there is a
constant  $N_{1} \in(0,\infty)$ such that
  $|v_{h}-g|\leq N_{1}\rho$
for   $h\in(0,h_{0}]$.

Now set $D=\Omega_{\kappa h}$. Then $D_{h}=\Omega_{(\kappa+1) h}$
and
we infer from Theorem \ref{theorem 6.5.1} that
for $|p|\leq h$
$$
|v_{h}(x+ l p)-v _{h}(x)|\leq N(|p|+
 \sup_{ \Omega^{c}_{(\kappa+1) h}}|v_{h}(\cdot+lp)-
v_{h}|)
$$
\begin{equation}
                                                \label{7.21.4}
\leq N( |p|+
  \sup_{ \Omega^{c}_{(\kappa+2) h}}|v_{h}-g|+
\sup_{ \Omega^{c}_{(\kappa+1) h}}|v_{h}-g|
+\sup_{ \Omega^{c}_{(\kappa+1) h}}|g(\cdot+lp)-
g|).
\end{equation}

Since
  $|v_{h}-g|\leq N_{1}\rho$,
$|v_{h}-g|\leq  N_{1}(\kappa+2) h$ outside $\Omega_{(\kappa+2)h}$,
and this along with \eqref{7.21.4} and the arbitrariness
of $l$ proves \eqref{3.26.6} for $|x-y|\leq h\leq h_{0}$,
which
finishes proving the corollary.

\mysection{Two results from \protect\cite{Kr12.1} and 
\protect\cite{Kr12.2}}
                                          \label{section 7.4.1}

We suppose that the assumptions in Section
\ref{section 1.22.1} are satisfied and
take ``cut-off'' functions 
$$
\eta\in C^{2}_{b}(\bR^{d}),\quad |\eta|\leq1,
\quad \zeta=\eta^{2}.
$$
 Fix an $h\in(0,\delta/(2K)]$ and 
set
$$
\Lambda_{h,1}=h\Lambda,\quad \Lambda_{h,n+1}= 
 \Lambda_{h,n}+h\Lambda ,\quad n\geq1,
\quad \Lambda_{h,\infty}=\bigcup_{n}\Lambda_{h,n}\,,
$$
 
Define
$$
Q^{o}=\{ x\in\Lambda_{h,\infty}:
x+3hB 
\subset \Omega\}=\Lambda_{h,\infty}\cap\Omega_{3h},
$$
$$
Q=\{x+\Lambda_{h,2}:x\in Q^{o}\},
\quad \delta   Q=Q\setminus Q^{o}   .
$$
 
Observe that $Q\subset\Omega_{h}$ and $Q$ is a finite set.
 The latter is due
to Assumption \ref{assumption 12.17.1}(iii) and follows from 
the fact that the number of points with integral coordinates
lying in a  bounded domain
is finite combined with the fact that there is a number
$M$ such that the coordinates of all points in 
$M\Lambda_{1,\infty}$
are integers. 
 
The following is a specification of Theorem  1.2
of \cite{Kr12.1} in the present setting.

\begin{theorem}
                                         \label{theorem 8.1.1}
Let $v$ satisfy \eqref{10.26.1}
in $Q $.
Then there is a constant $N_{1}\geq1$
depending only on $ d_{1},K$,  and $\delta$
such that, for any $\nu$ satisfying
$$ 
\nu\geq N_{1} (\sup_{\bR^{d}}|D^{2}\eta|
+ \sup_{\bR^{d}}|D \eta|^{2}+1),
$$
we have 
\begin{equation}
                                           \label{12.27.07}
\max_{k,Q}\zeta|\delta_{h,k}v|\leq N_{1}(\sqrt{\nu}
+\frac{1}{\nu}   ) \max_{Q}|v|+\frac{N_{1}}{\sqrt{\nu}} 
+\frac{N_{1}}{\nu} +N_{1}\max_{k,\delta Q}\zeta|\delta_{h,k}v|.
\end{equation}

\end{theorem}

\begin{remark}
In \cite{Kr12.1} a more general statement  than 
Theorem \ref{theorem 8.1.1} is proved under the assumption
that $a^{\alpha},b^{\alpha},c^{\alpha},f^{\alpha}$
are in $C^{0,1}(\bR^{d})$ rather than in $C^{1,1}(\bR^{d})$.
Also $\delta Q$ in \eqref{12.27.07} is replaced with a 
``thinner'' set and Assumption  \ref{assumption 12.17.1}(ii) 
is not used.
\end{remark}

Now comes a version of Theorem 1.1 of \cite{Kr12.2}.

\begin{theorem}
                                         \label{theorem 8.1.2}
Let $v$ satisfy \eqref{10.26.1}
in $Q $.
Then there exists a  constants $N_{2}$  depending only on 
  $ K,d_{1}$,
and $\delta$   
such that if   a constant $\nu>0$ satisfies
$$
\nu\geq N_{2}  (\sup_{\bR^{d}}|D^{2}\eta|
+ \sup_{\bR^{d}}|D \eta|^{2}+1),
$$
then  
$$
\max_{Q,i,j}\zeta  
|\delta_{h,i}\delta_{h,j}v|\leq \max_{\delta Q,i,j}\zeta  
|\delta_{h,i}\delta_{h,j}v| 
$$
$$
+N_{2}(\nu^{1/2} +\sup_{\bR^{d}}|D \eta| )
\max_{Q,i } 
|\delta_{h,i} v| + N_{2}(1+ \max_{Q}|v|).
$$
\end{theorem}

\mysection{Proof of Theorem \protect\ref{theorem 7.22.1}}
                                          \label{section 7.5.1}

(i) For $x\in\bR^{d}$ and $u=(u',u'')$, where
$$
u'=(u'_{-d_{1}},...,u'_{-1},u'_{0},u'_{1},
...,u'_{d_{1}}),\quad u''=(u''_{\pm1},...,u''_{\pm d_{1}}),
$$
introduce
$$
H(x,u)=
\sup_{\alpha\in A} 
 [a_{k}^{\alpha}( x)u''_{k}
+b_{k}^{\alpha}( x)u'_{k} 
-c^{\alpha}(x)u'_{0}+f^{\alpha}( x) ].
$$
By using the inequality
\begin{equation}
                                                      \label{7.12.5}
 \sup_{\alpha}F^{\alpha}-\sup_{\alpha}G^{\alpha} \leq
\sup_{\alpha}(F^{\alpha}-G^{\alpha})
\end{equation}
we easily conclude that $H(x,u)$ is locally Lipschitz continuous with
respect to $(x,u)$. As such it is almost everywhere differentiable.
Its derivatives are limits of finite differences and by using 
\eqref{7.12.5} again we see that at all points of differentiability
of $H$ we have
$$
|H_{u''_{j}}|,|H_{u' _{i}}|\leq K,\quad j=\pm1,...,\pm d_{1},
i=0,\pm1,...,\pm d_{1},
\quad H_{u'_{0}}\leq-\delta.
$$
 Furthermore,
$$
\sup_{\alpha}F^{\alpha}-\sup_{\alpha}G^{\alpha}
\geq\inf_{\alpha}(F^{\alpha}-G^{\alpha}),
$$
which implies that $H_{u''_{j}}\geq\delta$.

We have checked that $H$ (independent of $p$) satisfies
Assumption \ref{assumption 6.11.1}. Now assertion (i)
follows from Theorem \ref{theorem 7.12.1} and the fact that
$H(x,u)$ is a (Lipschitz) continuous in $x$, the latter being
again a consequence of \eqref{7.12.5} and
Assumption \ref{assumption 12.17.2}. 
This proves assertion (i) and combined with 
  Lemma \ref{lemma 7.4.1} shows that there is a
constant  $N  \in(0,\infty)$ such that
  $|v_{h}-g|\leq N \rho$
for   $h\in(0,h_{0}]$, which is part assertion (ii).

(ii) In light of \eqref{7.12.5}
at all points of  differentiability of $H(x,u)$
with respect to $(x,u)$
$$
|H_{x_{i}}(x,u)|\leq N(K,d_{1})(1+|u|)
$$
for all $i$. Hence,
owing to   Corollary \ref{corollary 7.21.1} to prove
assertion (ii),
it suffices to prove \eqref{6.25.3}, in which the first
estimate is obtained above.

Note that if $x\in\Lambda_{h,\infty}$ and 
 $x\not\in Q^{o}$, then $\rho(x)\leq 3h$, and
$$
|v_{h}(x+hl_{k})-g(x+hl_{k})|,  |v_{h}(x)-g(x)|\leq Nh
$$
 for any $k$. Hence for $x\in\Lambda_{h,\infty}
\setminus Q^{o}$
we have
$$
|\delta_{h,k}v_{h}(x)|\leq h^{-1}|v_{h}(x+hl_{k})-v_{h}(x)|
$$
$$
\leq h^{-1}( |v_{h}(x+hl_{k})-g(x+hl_{k})|
+|v_{h}(x)-g(x)|)+N\leq N,
$$
where $N$ is independent of $h$ and $x$. In short
\begin{equation}
                                       \label{6.25.2}
|\delta_{h,k}v_{h} |\leq N
\end{equation}
on $\Lambda_{h,\infty}\setminus Q^{o}$,
 where $N$ is independent of $h$.

By Theorem \ref{theorem 8.1.1} with $\eta\equiv1$
for all sufficiently small $h$
$$
\max_{k,Q}|\delta_{h,k}v_{h}|\leq N(1+\max_{\bR^{d}}|v_{h}|
+\max_{k,\delta Q}|\delta_{h,k}v_{h}|),
$$
where $N$ is independent of $h$. Since $|v_{h}-g|\leq N\rho$
and \eqref{6.25.2} holds on 
$\Lambda_{h,\infty}\setminus Q^{o}$, we conclude that
\eqref{6.25.2} holds on $\Lambda_{h,\infty}$, 
 provided that $h$
is sufficiently small, where $N$ is independent of $h$.
Actually, \eqref{6.25.2} holds on $\bR^{d}$ just because
any $x\in\bR^{d}$ can be placed into an appropriate shift
of $\Lambda_{h,\infty}$ with the shift not affecting any
of the above constants $N$.

 Thus, we estimated the second quantity
in \eqref{6.25.3}.

To estimate the last one we use Theorem \ref{theorem 8.1.2}.
Once again it suffices to concentrate on points in $\Lambda_{h,\infty}$.
Introduce $2\rho_{0}=\max\rho$.
It is a standard fact that for any $r\in(0,\rho_{0}]$ there exists an $\eta^{(r)}
\in C^{\infty}_{0}(\Omega)$ such that
$$
\eta^{(r)}=1\quad\text{on}\quad \Omega_{2r},\quad
\eta^{(r)}=0\quad\text{outside}\quad \Omega_{r},
$$
$$
|\eta^{(r)}|\leq1,\quad
|D\eta^{(r)}|\leq N/r,\quad|D^{2}\eta^{(r)}|\leq N/r^{2},
$$
where the constant $N$ is independent of $r$.
Furthermore, since $Q^{o}=\Lambda_{h,\infty}\cap\Omega_{3h}$
and $Q\subset\Omega_{h}$ it holds that
$\delta Q\subset \Omega\setminus \Omega_{3h}$. By taking
$r\geq 3h$, $r\leq1$, and
applying
Theorem \ref{theorem 8.1.2} with $\eta=\eta^{(r)}$ (and
$\nu=N/r^{2}$, where $N$ is an appropriate constant
independent of $r$ and $h$), 
 for sufficiently small $h>0$ and any $i,j$ 
we obtain on $\Lambda_{h,\infty}$
that 
\begin{equation}
                                       \label{6.25.4}
(\eta^{(r)}(x))^{2}|\delta_{h,i}\delta_{h,j}v_{h}(x)|\leq N/r,
\end{equation}
where $N$ is independent of $x,r$, and $h$.

Now we use the arbitrariness of $x$ and $r$. Take an $x\in\Lambda_{h,\infty}$
such that $\rho(x)\geq6h$ and take $r=\rho(x)/2$ ($\geq3h$).
Then $\eta^{(r)}(x)=1$ and \eqref{6.25.4} yields
$\rho(x)|\delta_{h,i}\delta_{h,j}v_{h}(x)|\leq N$ implying that
\begin{equation}
                                       \label{6.25.5}
(\rho(x)-6h)|\delta_{h,i}\delta_{h,j}v_{h}(x)|\leq N
\end{equation}
whenever $x\in\Lambda_{h,\infty}$ and
 $\rho(x)\geq6h$.
 However, \eqref{6.25.5} is obvious
if $\rho(x)\leq6h$. Thus \eqref{6.25.5} holds on $\Lambda_{h,\infty}$
and, as it was mentioned, on $\bR^{d}$. The theorem is
proved.

\mysection{Proof of Theorem \protect\ref{theorem 7.22.3}}
                                           \label{section 7.5.2}

We need two additional auxiliary results. Estimate \eqref{7.4.1} is quite similar
to (4.6) of \cite{DK05}
the latter being a particular case of the former
which occurs when $N_{1}^{*}=0$ (in (4.6) of \cite{DK05}
 it is assumed that
$u\in C^{0,1}$).

\begin{lemma}
                                          \label{lemma 7.1.1}
Let vector $l\in B$ be on the first basis axis in $\bR^{d}$.
Let $\varepsilon,h>0$,
$\eta\in C^{\infty}_{0}(\bR^{d})$ be a
spherically symmetric function
with support in $B_{\varepsilon}=\{x:|x|<\varepsilon\}$
    and let
$u$ be a bounded Borel function on $\bR^{d}$. Define
$$
B'_{\varepsilon}=\{z\in\bR^{d-1}:|z|<\varepsilon\}.
$$
Assume that for all $x,y \in[- 2h -\varepsilon ,  2h +\varepsilon   ]\times B'_{\varepsilon}$  
  we have
\begin{equation}
                                             \label{7.1.2}
|u(x )-u(y)|\leq N_{1} |x-y|+N_{1}^{*}h,\quad
|\Delta_{h,l}u(x)|\leq N_{2} .
\end{equation}
Introduce $w=u*\eta$. Then 
\begin{equation}
                                             \label{7.2.7}
|D_{l}^{2}w(0)-\Delta_{h,l}w(0)|\leq N_{2}h^{2} 
\|D^{2}_{l}\eta\|_{L_{1}}
+ h^{4}(N_{1}\varepsilon+N_{1}^{*}h) \|D^{6}_{l}\eta\|_{L_{1}}  .
\end{equation}
Also for $x=\gamma l$, $|x|\leq h$, we have
\begin{equation}
                                             \label{7.4.1}
|D_{l}^{2}w(x)| \leq N_{2}\|\eta\|_{L_{1}} + 
 h^{2}(N_{1}\varepsilon+N_{1}^{*}h) \|D^{4}_{l}\eta\|_{L_{1}}  ,
\end{equation}
Finally,
\begin{equation}
                                             \label{7.2.8}
|D_{l} w(0)-\delta_{h,l}w(0)|\leq
N_{2}h\|\eta\|_{L_{1}}
+ h^{3}(N_{1}\varepsilon+N_{1}^{*}h) \|D^{4}_{l}\eta\|_{L_{1}}.
\end{equation}

\end{lemma}

Proof. First for $n=1,2,...,6$ and $x=\gamma l  $, $|x|\leq 2h$,
we have
$$
D^{n}_{l}w(x)=\int_{B_{\varepsilon}}u(x-y)D^{n}_{l}\eta(y)\,dy=
\int_{B_{\varepsilon}}[u(x-y)-u(x)]D^{n}_{l}\eta(y)\,dy.
$$
Hence,
\begin{equation}
                                             \label{7.1.3}
|D^{n}_{l}w(x)|\leq (N_{1}\varepsilon+N_{1}^{*}h)
\|D^{n}_{l}\eta\|_{L_{1}}.
\end{equation}

Next observe that 
for smooth functions $f(t)$ of one real variable $t$ we have
$$
f''(0)=\frac{f(h)-2f(0)+f(-h)}{h^{2}}-\frac{h^{2}}{6}
\int_{-1}^{1}(1-|t|)^{3}f^{(4)}(th)\,dt.
$$
Hence
\begin{equation}
                                             \label{7.1.4}
|D_{l}^{2}w(0)-\Delta_{h,l}w(0)|\leq  h^{2} 
\sup_{|t_{1}|\leq1}|D_{l}^{4}w(lt_{1}h)|.
\end{equation}
Therefore, to prove \eqref{7.2.7} it suffices to show that
for $|t_{1}|\leq1$
\begin{equation}
                                             \label{7.2.08}
|D_{l}^{4}w(lt_{1}h)|\leq N _{2} 
\|D ^{2}_{l}\eta\|_{L_{1}}
+ h^{2}(N_{1}\varepsilon+N_{1}^{*}h)\|D^{6}_{l}\eta\|_{L_{1}}.
\end{equation}
To this end we use \eqref{7.1.4} which implies that
\begin{equation}
                                             \label{7.2.1}
|D_{l}^{2} u*\eta (0)|\leq |\Delta_{h,l} u*\eta (0)|+  h^{2} 
\sup_{|t_{1}|\leq1}|D_{l}^{2} u*D_{l}^{2}\eta (lt_{1}h)|.
\end{equation}
By applying \eqref{7.2.1} to $D^{2 }_{l}\eta$ 
 in place of $\eta$ and $lt_{1}h$ in place of $0$
  we find
$$
|D_{l}^{4}w(lt_{1}h)|=|D_{l}^{2} u*D^{2 }_{l}\eta (lt_{1}h)|
\leq |\Delta_{h,l} u*D^{2 }_{l}\eta (lt_{1}h)|
$$
$$
+  h^{2} 
\sup_{|t_{2}|\leq1}|D_{l}^{2} u*D_{l}^{4}\eta (l(t_{1}+t_{2})h)|,
$$
Here, owing to \eqref{7.1.2} and the fact that $|t_{1}+t_{2}|h
+\varepsilon\leq2h+\varepsilon$,
 $$
|\Delta_{h,l} u*
D^{2}_{l}\eta (l(t_{1}+t_{2})h)|\leq N_{2}\|D^{2}_{l}\eta\|_{L_{1}}.
$$
Furthermore, by \eqref{7.1.3} and the fact that $|t_{1}+t_{2}|h
 \leq2h $
$$
|D_{l}^{2} u*D_{l}^{4}\eta (l(t_{1}+t_{2} )h)|=
|  u*D_{l}^{6}\eta (l(t_{1}+t_{2} )h)|
$$
$$
 \leq
(N_{1}\varepsilon+N_{1}^{*}h)\|D^{6}_{l}\eta\|_{L_{1}}.
$$
This proves \eqref{7.2.08} and \eqref{7.2.7}.

Since, as above, $|\Delta_{h,l}w(0)|\leq N_{2}\|\eta\|_{L_{1}}$,
we obtain \eqref{7.4.1} for $x=0$ directly from
\eqref{7.1.4} and \eqref{7.1.3}. For other values of $x$
we obtain \eqref{7.4.1} upon observing that the above argument
is valid if we replace $0$ and $lt_{1}h$ in
\eqref{7.1.4} with $x$ and $x+lt_{1}h$,  respectively.

To prove \eqref{7.2.8}
observe that in the 
one-dimensional case
$$
f'(0)=\frac{f(h)-f(0)}{ h}-h
\int_{0}^{1}(1- t ) f''(th)\,dt.
$$
It follows that
\begin{equation}
                                                        \label{7.8.6}
|D_{l}w(0)- \delta_{h,l} w(0)|\leq
h
\int_{0}^{1}(1- t ) |D_{l}^{2}w (lth)|\,dt
\leq h\sup_{t\in[0,1]}|D_{l}^{2}w (lth)|,
\end{equation}
which along with \eqref{7.4.1} yields \eqref{7.2.8}.
The lemma is proved.  

\begin{remark}
                                         \label{remark 7.1.1}
Take a nonnegative spherically symmetric 
$\zeta\in C^{\infty}_{0}(\bR^{d})$
whose integral is one and whose support is in $B$
 and for $\varepsilon>0$ and
locally integrable functions $u$ on $\bR^{d}$ use the notation
$$
u^{(\varepsilon)}(x)=\int_{\bR^{d}}u(x-\varepsilon y)\zeta(y)\,dy.
$$
 
Denote by $v_{h} $
a unique bounded solution 
of \eqref{10.26.1}, with zero boundary data on $\partial_{h}\Omega$.
By Theorem \ref{theorem 7.22.1} the function $v_{h} $
is well defined at least for sufficiently small $h>0$.  

We want to explain why we iterated \eqref{7.2.1} and why
only once. We will apply Lemma \ref{lemma 7.1.1} to $u=v_{h}$. 
For simplicity assume that $b\equiv0$ and notice that on $\Omega_{h}$
for any $\alpha\in A$ we have
\begin{equation}
                                                   \label{7.15.1}
  a^{\alpha}_{k} D^{2}_{l_{k}}v  
-c^{\alpha} v +f^{\alpha} \leq0.
\end{equation}
Next assume that $a $  and $c$ are independent of $x$
as in \cite{DK05}.
Then we can mollify all terms in \eqref{7.15.1} and obtain that
in $\Omega_{h+\varepsilon}$ we have
\begin{equation}
                                                   \label{7.15.2}
a_{k}^{\alpha} \Delta_{h,k}v _{h}^{ (\varepsilon)} 
-c^{\alpha} v _{h}^{ (\varepsilon)} +[f^{\alpha} ]
^{ (\varepsilon)}\leq0.
\end{equation}
It is well known  that
$|[f^{\alpha} ]
^{ (\varepsilon)}-f^{\alpha}|\leq N\varepsilon^{2}$ since 
$f^{\alpha}\in C^{1,1}$. Therefore, \eqref{7.15.2} implies that
$$
a_{k}^{\alpha} \Delta_{h,k}v _{h}^{ (\varepsilon)} 
-c^{\alpha} v _{h}^{ (\varepsilon)} + f^{\alpha} 
 \leq N\varepsilon^{2}.
$$
Now we replace $\Delta_{h,k}v _{h}^{ (\varepsilon)}$ 
 with $D^{2}_{l_{k}}
v _{h}^{ (\varepsilon)}$. Ignoring for a moment that our estimates of
$\Delta_{h,k}v_{h}$ are local in $\Omega$, we obtain from
Lemma \ref{lemma 7.1.1} that
$$
a_{k}^{\alpha} D^{2}_{l_{k}}v _{h}^{ (\varepsilon)} 
-c^{\alpha} v _{h}^{ (\varepsilon)} + f^{\alpha} 
 \leq N[\varepsilon^{2}+ h^{2}\varepsilon^{-2}+h^{4}(\varepsilon+h)
\varepsilon^{-6}]=:NM(h,\varepsilon).
$$

It follows from Lemma \ref{lemma 7.7.1}, that
for, perhaps, another constant $N$
$$
H[v _{h}^{ (\varepsilon)}+N \Psi M(h,\varepsilon)]\leq 0
$$
in $\Omega_{h+\varepsilon}$
By the maximum principle
\begin{equation}
                                                  \label{7.17.1}
v\leq v _{h}^{ (\varepsilon)}+N M(h,\varepsilon)
+\sup_{\Omega_{h+\varepsilon}^{c}}(v-v _{h}^{ (\varepsilon)})_{+}.
\end{equation}
It seems that the best possible estimate of the last term is that
it is less than $N(h+\varepsilon)$ and we will prove this estimate
later. Then
$$
v\leq v _{h}^{ (\varepsilon)}+N [M(h,\varepsilon)+h+\varepsilon],
$$
\begin{equation}
                                                           \label{7.16.1}
v\leq v_{h}+N [M(h,\varepsilon)+h+\varepsilon]
+\sup(v _{h}^{ (\varepsilon)}
-v_{h})_{+}.
\end{equation}
Let us ignore the contribution of the last term in the right-hand side
and try to make $M(h,\varepsilon)+h+\varepsilon$ as mall as
possible on the  account of arbitrariness in choosing
$\varepsilon$. Observe that this quantity contains
$\varepsilon+h^{2}\varepsilon^{-2}$ which is bigger than $\gamma
h^{2/3}$, where the constant $\gamma>0$ is independent of $h$.
Furthermore, $\varepsilon+h^{2}\varepsilon^{-2}=2h^{2/3}$
when $\varepsilon=h^{2/3}$. With this $\varepsilon$
$$
h^{4}(\varepsilon+h)
\varepsilon^{-6}=O(h^{2/3})
$$
as well and we obtain that $v\leq v_{h}+Nh^{2/3}$ for sufficiently
small $h>0$.

This is roughly the way we are going to use \eqref{7.2.7}.

Now imagine that we did not enhance the estimate
of $D^{4}_{l}w$ and  after \eqref{7.1.4}
just used \eqref{7.1.3} to obtain
$$
|D_{l}^{2}v _{h}^{ (\varepsilon)}-\Delta_{h,l}v _{h}^{ (\varepsilon)}|
\leq  Nh^{2}(\varepsilon+h)\varepsilon^{-4}.
$$
Then we would obtain \eqref{7.16.1} with
$$
M(h,\varepsilon)=\varepsilon^{2}+h^{2}(\varepsilon+h)
\varepsilon^{-4}$$
and $M(h,\varepsilon)+h+\varepsilon$ would contain the term
$\varepsilon+h^{2}\varepsilon^{-3}$ whose minimum with respect
to $\varepsilon$ is of order $h^{1/2}$
and is comparable with its value at $\varepsilon=h^{1/2}$.
Then at best we would have that $v\leq v_{h}+Nh^{1/2}$.

On the other hand, we could iterate \eqref{7.2.1} one more time
and obtain that
$$
|D_{l}^{2} u*D_{l}^{4}\eta (l(t_{1}+t_{2})h)|\leq
|\Delta_{h,l} u*D_{l}^{4}\eta (l(t_{1}+t_{2})h)|
$$
$$
+h^{2}\sup_{|t_{3}|\leq1}|D_{l}^{2} u*D_{l}^{6}\eta (l(t_{1}+t_{2}+t_{3})h)|,
$$
which in our case  leads to
$$
|D_{l}^{2}v _{h}^{ (\varepsilon)}-\Delta_{h,l}v _{h}^{ (\varepsilon)}|
\leq Nh^{2}[\varepsilon^{-2}+h^{2}(\varepsilon^{-4}+
h^{2}(\varepsilon+h)\varepsilon^{-8})]
$$
$$
=N(h^{2}\varepsilon^{-2}+h^{4}\varepsilon^{-4}+h^{6}
(\varepsilon+h)\varepsilon^{-8}).
$$
This time again $M(h,\varepsilon)+h+\varepsilon $  contains 
$\varepsilon +h^{2}\varepsilon^{-2}$, which
  will not lead to a better rate than $h^{2/3}$.

It is seen from the above that the fact that 
$|[f^{\alpha} ]
^{ (\varepsilon)}-f^{\alpha}|\leq N\varepsilon^{2}$ plays no significant role.
The estimate $|[f^{\alpha} ]
^{ (\varepsilon)}-f^{\alpha}|\leq N\varepsilon $ would do equally well.
Also in this framework we will be satisfied with
estimating the last term in the right-hand side of \eqref{7.16.1} just
by $N\varepsilon$, which is quite easy.

It is worth noting that the situation with constant in $x$
coefficients $a,b,$ and $c$ is quite different in the whole
space (see \cite{DK05}).
 There no boundary term like the last term in \eqref{7.17.1}
appears and one ends up with $v\leq v_{h}^{(\varepsilon)}+NM(h,\varepsilon)$.
With some additional effort (even in our setting of bounded smooth
$\Omega$ and variable coefficients, see
Lemma \ref{lemma 7.2.1}) one can prove that 
$|v_{h}^{(\varepsilon)}-v_{h}|\leq N(\varepsilon^{2}+h)$.
Then
\begin{equation}
                                                       \label{7.17.2}
v\leq v_{h} +N
[h+\varepsilon^{2}+ h^{2}\varepsilon^{-2}+h^{4}(\varepsilon+h)
\varepsilon^{-6}].
\end{equation}
The minimum of $\varepsilon^{2}+ h^{2}\varepsilon^{-2}$ with respect
to $\varepsilon$ is proportional to its value at $\varepsilon=h^{1/2}$
and is of order $h$. Other error terms on the right in \eqref{7.17.2}
are of the same or higher order. In this way it is proved in 
\cite{DK05} that $v\leq v_{h} +Nh$. It is also shown there that
in general this estimate is optimal, so that there is no need
to even consider additional iterations of \eqref{7.2.1}.

Finally, we point out that even for the equations in the whole
space with {\em variable\/} coefficients the error term
of order $\varepsilon$ still appears in the transition from
\eqref{7.15.1} to \eqref{7.15.2}. The method of
``shaking the coefficients'' produces an error of the same
order, which allows us not to use this method on the
account that we have a good control of the second-order differences 
of $v_{h}$.

 \end{remark}

\begin{lemma}
                                           \label{lemma 6.29.1}
There exists a constant $N$ such that for all sufficiently
small $h>0$, for any $\varepsilon>0$,

(i) In $\bR^{d}$ we have $|v^{(\varepsilon)}_{h}-v_{h}|
\leq N( \varepsilon +h)$;

(ii) In $\Omega_{4\varepsilon+16h}$ for any $\alpha\in A$ we have
$$
  a_{k}^{\alpha} D^{2}_{l_{k}}v^{ (\varepsilon)}_{h} 
+b_{k}^{\alpha}  D _{l_{k}}v^{ (\varepsilon)}_{h} 
-c^{\alpha} v^{ (\varepsilon)}_{h}(x)+f^{\alpha}   
$$
\begin{equation}
                                                    \label{6.29.2}
\leq
N(h+\varepsilon+h^{2}\varepsilon^{-2})\rho^{-1}+N
h^{3}\varepsilon^{-6}(\varepsilon+h)(\varepsilon^{2}+h).
\end{equation}

\end{lemma}

Proof. Assertion (i) follows immediately from the fact that
$|v_{h}(x)-v_{h}(y)|\leq N(|x-y|+h)$ (see
Theorem \ref{theorem 7.22.1}).

(ii) Fix an $\alpha\in A$ and
observe that for $x\in\Omega_{h+\varepsilon}$ we have
\begin{equation}
                                                    \label{7.4.5}
[ a_{k}^{\alpha} \Delta_{h,k}v _{h} 
+b_{k}^{\alpha} \delta_{h,k}v _{h} 
-c^{\alpha} v _{h} +f^{\alpha} ]
^{ (\varepsilon)} (x)\leq0.
\end{equation}

Next,
$$
[ a_{k}^{\alpha} \Delta_{h,k}v _{h} ]
^{ (\varepsilon)} (x)=a_{k}^{\alpha}(x)
\Delta_{h,k}v^{ (\varepsilon)}_{h}  
 (x)
$$
$$
+\int_{B }[a_{k}^{\alpha}(x
+\varepsilon y)-a_{k}^{\alpha}(x)]\zeta(y)
\Delta_{h,k}v _{h}(x
+\varepsilon y)\,dy,
$$
where owing to Theorem \ref{theorem 7.22.1},
for $x\in\Omega_{6h+\varepsilon}$ the last term by magnitude
is less than
$$
N\varepsilon\int_{B }\frac{1}{\rho(x+\varepsilon y)-6h}\,dy,
$$
which is less than $N\varepsilon \rho^{-1}(x)$ if
$x\in\Omega_{12h+2\varepsilon}$ since then
$$
\rho(x+\varepsilon y)-\rho(x)\geq-\varepsilon,
\quad  \rho(x+\varepsilon y)-(1/2)\rho(x)\geq(1/2)\rho(x)
-\varepsilon\geq 6h,
$$
$$
\rho(x+\varepsilon y)-6h\geq(1/2)\rho(x).
$$
Hence in $\Omega_{12h+2\varepsilon}$
$$
a_{k}^{\alpha} 
\Delta_{h,k}v^{ (\varepsilon)}_{h}  \leq 
[ a_{k}^{\alpha} \Delta_{h,k}v _{h} ]
^{ (\varepsilon)} + N\varepsilon \rho^{-1}.
$$

Then again owing to Theorem \ref{theorem 7.22.1}
in $\bR^{d}$
$$
b_{k}^{\alpha} 
\delta_{h,k}v^{ (\varepsilon)}_{h} \leq 
[ b_{k}^{\alpha} \delta_{h,k}v _{h} ]
^{ (\varepsilon)}  + N\varepsilon,
\quad -c^{\alpha} 
 v^{ (\varepsilon)}_{h} \leq -[ c^{\alpha} v _{h} ]
^{ (\varepsilon)}   + N\varepsilon,
$$
$$
 f^{\alpha}  \leq [ f^{\alpha} ]
^{ (\varepsilon)} 
 + N\varepsilon .
$$
Coming back to \eqref{7.4.5} we conclude that in $\Omega_{12h+2\varepsilon}$
we have
\begin{equation}
                                                    \label{7.4.8}
  a_{k}^{\alpha} \Delta_{h,k}v _{h}^{ (\varepsilon)} 
+b_{k}^{\alpha} \delta_{h,k}v _{h}^{ (\varepsilon)} 
-c^{\alpha} v _{h}^{ (\varepsilon)} +f^{\alpha}  
 (x)\leq N\varepsilon \rho^{-1}.
\end{equation}

Now we are going to use Lemma \ref{lemma 7.1.1} in order to replace
$\Delta_{h,k}$ and $\delta_{h,k}$ in \eqref{7.4.8} with
$D^{2}_{l_{k}}$ and  $D_{l_{k}}$, respectively. Of course,
this time we take   $\eta(x)=
\varepsilon^{-d}\zeta(x/\varepsilon)$ and $u=v_{h}$
in Lemma \ref{lemma 7.1.1}. This lemma is stated only for 
vectors on the first basis axis.  The reader understands that
one can prove an appropriately modified
statements for any vector $l\in B$ 
under the assumption that 
\eqref{7.1.2} holds in
the cylinder $C_{2\varepsilon+4h,\varepsilon,l}$ centered at the origin
with hight $2\varepsilon+4h$, base radius $\varepsilon$,
and axis parallel to $l$. Naturally, one can move
such cylinders to be centered at any point $x_{0}$.
Observe that if $x_{0}\in\Omega_{2\varepsilon+8h}$,
then, for any $l\in\Lambda$, all points of the   cylinder
$x_{0}+C_{2\varepsilon+4h,\varepsilon,l}$
are at a distance not less than 
$$
\rho(x_{0})-2\varepsilon-2h>6h
$$
 from $\Omega^{c}$.
Hence, while applying Lemma \ref{lemma 7.1.1}
to $x_{0}+C_{2\varepsilon+4h,\varepsilon,l}$ we can take
the constant $N_{2}$ to be 
$$
N(\rho(x_{0})-2\varepsilon-8h)^{-1},
$$ 
where $N$ is the constant from Theorem \ref{theorem 7.22.1}.
Notice that 
$$
\rho(x_{0})-2\varepsilon-8h\geq{1/2}\rho(x_{0})
$$
if $x_{0}\in \Omega_{4\varepsilon+16h}$. Therefore,
on $\Omega_{4\varepsilon+16h}$
by Lemma \ref{lemma 7.1.1} we obtain
$$
|D_{l_{k}}v_{h}^{(\varepsilon)}-\delta_{h,k}
v_{h}^{(\varepsilon)}|\leq Nh\rho^{-1}+Nh^{3}(\varepsilon+h)
\varepsilon^{-4},
$$
$$
|D^{2}_{l_{k}}v_{h}^{(\varepsilon)}-\Delta_{h,k}
v_{h}^{(\varepsilon)}|\leq Nh^{2}\varepsilon^{-2}\rho^{-1}
+Nh^{4}(\varepsilon+h)
\varepsilon^{-6}.
$$
These estimates and \eqref{7.4.8} yield
\eqref{6.29.2}. The lemma is proved.

Assertion (i) of Lemma \ref{lemma 6.29.1}
that $|v^{(\varepsilon)}_{h}-v_{h}|
\leq N( \varepsilon +h)$ in $\bR^{d}$ can be
improved for points at a fixed distance
from $\Omega^{c}$ if one uses the following lemma, in which
  $e_{1},...,e_{d}$ is the standard orthonormal
basis in $\bR^{d}$. 

\begin{lemma}
                                        \label{lemma 7.2.1}
Let
$u$ be a bounded Borel function on $\bR^{d}$.
Assume that for all $x,y \in[-(2h+\varepsilon ), 
2h+\varepsilon ]^{d}$  
  we have
$$
|u(x )-u(y)|\leq N_{1}(|x-y|+h),
$$
$$
|\Delta_{h,e_{i}}u(x)|\leq N_{2}\quad i=1,...,d. %%%
$$
Then for $\varepsilon\geq h$
 \begin{equation}
                                                      \label{7.4.2}
|u^{(\varepsilon)}(0)-u(0)|\leq N(\varepsilon^{2}+h),
\end{equation}
where  the constant  $N$ 
depends only on $\zeta,d,N_{1}$, and $N_{2}$.
\end{lemma}
 
Proof. We follow closely the argument in Remark 4.5
of \cite{DK05}. For $r>0$ introduce
$$
w_{r}(x)=r^{-d/2}\int_{\bR^{d}}\zeta_{r}(x-y)u(y)\,dy,
$$
where
$$
\zeta_{r}(x)=\zeta(xr^{-1/2}).
$$
Simple computations show that
$$
\frac{\partial}{\partial r}\big[\frac{1}{r^{d/2}}\zeta_{r}(x)\big]=
\frac{1}{r^{d/2}}D^{2}_{i}\zeta_{i,r}(x),
$$
where (with no summation in $i$)
$$
\zeta_{i,r}(x)=\zeta_{i}(xr^{-1/2}),\quad
\zeta_{i}(x)=-\frac{1}{2}\int_{-\infty}^{x_{i}}
\zeta(x-x_{i}e_{i}+se_{i})s\,ds.
$$
Observe that since $\zeta$ is spherically symmetric,
the support of $\zeta_{i}$ is in $B $ and, of course,
$\zeta_{i}\in C^{\infty}_{0}(\bR^{d})$. Also notice for the future that
(no summation in~$i$)
$$
|D^{4}_{i}\zeta_{i,r}(x)|=\frac{1}{r^{2}}|D^{4}_{i}\zeta_{i}|
 (xr^{-1/2} )\leq \frac{1}{r^{2}}(
 |D^{3}_{i}\zeta|+|D^{2}_{i}\zeta|)(xr^{-1/2} ).
$$

It follows that  
$$
w_{t}(0)-w_{r}(0) =\int_{r}^{t}
\frac{1}{s^{d/2}}D^{2}_{i}[u*\zeta_{i,s}](0)\,ds.
$$
The support of $\zeta_{i,s}$ lies in $B_{\sqrt{s}}$.
Therefore, if    $0<r<t\leq\varepsilon^{2}$ 
and $s\in[r,t]$, we can use \eqref{7.4.1} 
and obtain (no summation in $i$)
$$
|D^{2}_{i}[u*\zeta_{i,s}](0)|\leq N_{2}\|\zeta_{i,s}\|_{L_{1}}
+N_{1}
 h^{2}(s^{1/2}+h) \|D^{4}\zeta_{i,s}\|_{L_{1}}
$$
$$
=N_{2}s^{d/2}\|\zeta_{i }\|_{L_{1}}
+N_{1}
 h^{2}(s^{1/2}+h)s^{d/2-2} \|D^{4}\zeta_{i }\|_{L_{1}}.
$$

Hence
$$
|w_{t}(0)-w_{r}(0)|\leq N\int_{r}^{t}(1+h^{2}s^{-3/2}+h^{3}s^{-2})
\,ds
$$
$$
\leq N(t-r)+Nh^{2}r^{-1/2}+Nh^{3}r^{-1},
$$
where and below the constants $N$ 
depend only on $\zeta,d,N_{1}$, and $N_{2}$. We combine this
with
$$
|w_{r}(0)-u(0)|=\big|\int_{\bR^{d}}\zeta(y)[u( yr^{1/2})-u(0)]
\,dy\big|\leq N_{1}(r^{1/2}+h).
$$
Then we get
$$
|w_{t}(0)-u(0)|
\leq N(t-r)+Nh^{2}r^{-1/2}+Nh^{3}r^{-1}+N(r^{1/2}+h),
$$
 which leads to \eqref{7.4.2} if we take $t=\varepsilon^{2}$
and $r=h^{2}$. The lemma is proved.

{\bf Proof of Theorem  \ref{theorem 7.22.3}}. Take
a constant $\mu\geq 1$ from Lemma \ref{lemma 7.7.1}.
By Lemmas \ref{lemma 6.29.1}(ii) and \ref{lemma 7.7.1}, 
there is a constant
$N$ independent of $\alpha$, $h$, and $\varepsilon$ 
such that for 
$$
N_{1}:= N(h+\varepsilon+h^{2}\varepsilon^{-2}) +N
h^{3}\varepsilon^{-6}(\varepsilon+h)(\varepsilon^{2}+h),
$$
$$
w^{\varepsilon}_{h}:=v^{(\varepsilon)}_{h}+N_{1}\Phi
$$
and sufficiently small $h>0$
we have
$$
 a_{k}^{\alpha}  D^{2}_{l_{k}}w^{\varepsilon}_{h}
+b_{k}^{\alpha} D _{l_{k}}w^{\varepsilon}_{h} 
-c^{\alpha}w^{\varepsilon}_{h}+f^{\alpha}  \leq0
$$
for all $\alpha\in A$ in $\Omega_{\kappa}$, where $\kappa=
(\mu h)\vee(4\varepsilon+16 h)$. By the maximum principle
\begin{equation}
                                                        \label{7.8.1}
v \leq w^{\varepsilon}_{h}
+\max_{\Omega\setminus\Omega_{\kappa}}(v-w^{\varepsilon}_{h})_{+}
\leq v_{h}^{(\varepsilon)}+N_{1}\Phi+
\max_{\Omega\setminus\Omega_{\kappa}}(v-v_{h}^{(\varepsilon)})_{+}
\end{equation}
in $\Omega$. By Theorem  \ref{theorem 7.22.2} we have
$|v-g|\leq N\kappa$ in $\Omega\setminus\Omega_{\kappa}$.
Furthermore, by  Lemma 
\ref{lemma 6.29.1}(i) we have $|v^{(\varepsilon)}_{h}-v_{h}|
\leq N( \varepsilon +h)$ everywhere which
implies that 
$$
|v^{(\varepsilon)}_{h}-g|\leq |v_{h}-g|
+N( \varepsilon +h)
$$
 and along with 
Theorem \ref{theorem 7.22.1} yields that
$$
|v^{(\varepsilon)}_{h}-g|\leq N\kappa 
+N( \varepsilon +h)
$$
 in $\Omega\setminus\Omega_{\kappa}$.
Upon combining this with \eqref{7.8.1} and observing that
$\kappa\leq N(\varepsilon+h)$ we get
$$
v  
\leq v_{h} + 
N(h+\varepsilon+h^{2}\varepsilon^{-2}) +N
h^{3}\varepsilon^{-6}(\varepsilon+h)(\varepsilon^{2}+h),
$$
 which yields that in $\Omega$ for sufficiently small $h>0$ we have
\begin{equation}
                                                        \label{7.8.03}
v  
\leq v_{h} + Nh^{2/3}
\end{equation}
if we set $\varepsilon=h^{2/3}$.

To prove that 
\begin{equation}
                                                        \label{7.9.1}
v_{h}\leq v+Nh^{2/3}
\end{equation}
 we reverse the roles
of $v_{h}$ and $v$. In $\Omega_{\varepsilon}$ we have
\begin{equation}
                                           \label{3.13.1}
  [a_{k}^{\alpha} D^{2}_{l_{k}}v  
+b_{k}^{\alpha}  D _{l_{k}}v 
-c^{\alpha} v +f^{\alpha}]^{(\varepsilon)}\leq0.
\end{equation}
Furthermore, for $x\in\Omega_{\varepsilon}$
$$
a^{\alpha}_{k}D^{2}_{l_{k}}v^{(\varepsilon)}(x)
\leq [a_{k}^{\alpha} D^{2}_{l_{k}}v]^{(\varepsilon)}
+N\varepsilon\sup_{k,x+\varepsilon B} |D^{2}_{l_{k}}v|,
$$
where by Theorem \ref{theorem 7.22.2} the second term on  the right 
is dominated by $N(\rho(x)-\varepsilon)^{-1}$, which is less
than $N\rho^{-1}(x)$ if $\rho(x)-\varepsilon\geq(1/2)\rho(x)$,
 that is if $x\in\Omega_{2\varepsilon}$.
Similarly one estimates $ b_{k}^{\alpha}  D _{l_{k}}v^{(\varepsilon)} $,
 $-c^{\alpha}v^{(\varepsilon)} $ and $[f^{\alpha}]^{(\varepsilon)}$.
 Then one concludes from \eqref{3.13.1}
that
\begin{equation}
                                                        \label{7.8.4}
  a_{k}^{\alpha} D^{2}_{l_{k}}v^{(\varepsilon)}  
+b_{k}^{\alpha}  D _{l_{k}}v^{(\varepsilon)} 
-c^{\alpha}v^{(\varepsilon)} + f^{\alpha} \leq
N\varepsilon\rho^{-1} 
\end{equation}
in $\Omega_{2\varepsilon}$.

Next, as in \eqref{7.1.4}
for  $x\in\Omega_{2h+2\varepsilon}$ we have
$$
|D^{2}_{l_{k}}v^{(\varepsilon)}(x)-
\Delta_{h,k}v^{(\varepsilon)}(x)|\leq Nh^{2}
\sup_{x+hB}|D_{l_{k}}^{4}v^{(\varepsilon)}(y)|,
$$
where 
$$
|D_{l_{k}}^{4}v^{(\varepsilon)}(y)|=
|D_{l_{k}}^{2}[D_{l_{k}}^{2}v]^{(\varepsilon)}(y)|
\leq N\varepsilon^{-2}(\rho(x)-h-\varepsilon)^{-1}
\leq N\varepsilon^{-2} \rho^{-1}(x)
$$
if $y\in x+hB$.
Also as in \eqref{7.8.6}
$$
|D _{l_{k}}v^{(\varepsilon)}(x)-
\delta_{h,k}v^{(\varepsilon)}(x)|\leq Nh 
\sup_{x+hB}|D_{l_{k}}^{2}v^{(\varepsilon)}|,
$$
where 
$$
|D_{l_{k}}^{2}v^{(\varepsilon)}(y)| 
\leq N (\rho(x)-h-\varepsilon)^{-1}
\leq N \rho^{-1}(x) 
$$
if $y\in x+hB$ and $x\in\Omega_{2h+2\varepsilon}$.
Hence \eqref{7.8.4} implies that 
$$
 a_{k}^{\alpha} \Delta_{h,k}v^{(\varepsilon)}  
+b_{k}^{\alpha} \delta_{h,k}v^{(\varepsilon)} 
-c^{\alpha}v^{(\varepsilon)} + f^{\alpha} \leq
N(\varepsilon + h^{2}\varepsilon^{-2}  
+ h)\rho^{-1}
$$
in $\Omega_{2h+2\varepsilon}$. At this point it is convenient
to extend $v$ outside $\Omega$ as $g$. By using
  Lemma \ref{lemma 7.7.1} and the maximum principle as above we obtain that
in $\bR^{d}$
$$
v_{h}\leq v^{(\varepsilon)}+
N(\varepsilon + h^{2}\varepsilon^{-2}  
+ h)+\sup_{\partial_{\chi} \Omega}(v_{h}-v^{(\varepsilon)})_{+},
$$
where $\chi=(\mu h)\vee(2h+2\varepsilon)$. Furthermore, in
$\Omega\setminus \Omega_{\chi}$ we have $|v_{h}-g|\leq N\chi$, which
follows from Theorem \ref{theorem 7.22.1}(ii) and
$|v^{(\varepsilon)}-g|\leq N(\chi+\varepsilon)\leq N\chi$,
which follows from the fact that $v\in C^{0,1}(\bar{\Omega})$
by Theorem \ref{theorem 7.22.2}. Since $\chi\leq N(h+\varepsilon)$
we conclude that
$$
v_{h}\leq v^{(\varepsilon)}+
N(\varepsilon + h^{2}\varepsilon^{-2}  
+ h).
$$
The boundedness of the first derivatives of $v$ implies that
$v^{(\varepsilon)}\leq v+N\varepsilon$. Hence,
$$
v_{h}\leq v +
N(\varepsilon + h^{2}\varepsilon^{-2}  
+ h),
$$
which yields \eqref{7.9.1} for $\varepsilon=h^{2/3}$
and along with \eqref{7.8.03} brings the proof
of the theorem to an end.


\begin{thebibliography}{mm}



\bibitem{BJ02} G. Barles and E.R. Jakobsen,  {\em
  On the convergence rate of
 approximation schemes for Hamilton-Jacobi-Bellman equations\/},
 M2AN Math. Model. Numer. Anal., Vol. 36 (2002), 
No. 1, 33--54.

\bibitem{BJ05} G. Barles and E.R. Jakobsen,  {\em Error bounds for monotone 
approximation schemes for Hamilton-Jacobi-Bellman equations\/},
 SIAM J. Numer. Anal., Vol. 43 (2005), No. 2, 
540--558.

\bibitem{BJ07} G. Barles and E.R. Jakobsen, 
{\em Error bounds for monotone
approximation schemes for parabolic Hamilton-Jacobi-Bellman 
equations\/}. Math.
Comp., Vol. 76 (2007), No. 260, 1861--1893.

\bibitem{BS} G. Barles and P.E. Souganidis, {\em
Convergence of approximation 
schemes for fully nonlinear second
order equations\/}, Asymptotic Anal., Vol. 4 (1991), No. 3, 271--283.


\bibitem{CL} M.G. Crandall and P.-L. Lions,  {\em
Two approximations of solutions of Hamilton-Jacobi equations\/},
Math. Comp., Vol. 43 (1984), No. 167, 1--19.


\bibitem{DK05} Hongjie Dong and N.V. Krylov,
 {\em
On the rate of convergence
of finite-difference approximations for Bellman equations with
constant coefficients\/},  Algebra i Analiz, Vol. 17 (2005),
 No. 2, 108-132; St.~Petersburg
Math. J., Vol. 17 (2006), No. 2, 295--313.


\bibitem{DK07} Hongjie Dong and N.V. Krylov,
{\em  The rate of convergence
of finite-difference approximations for parabolic Bellman equations
with Lipschitz coefficients in cylindrical
domains\/},   Applied Math. and Optimization, Vol. 56 (2007), No. 1,
37--66.

\bibitem{FS06} W.L. Fleming and H.M. Soner,  
 ``Controlled Markov processes and
viscosity solutions", second edition, Stochastic Modelling and
Applied Probability, Vol. 25, Springer, New York, 2006.

\bibitem{Kr77} N.V. Krylov,  ``Controlled diffusion processes'', 
Nauka, Moscow,  1977 in Russian; English translation:
 Springer,
1980.

\bibitem{Kr72} N.V. Krylov, {\em
   On control of the solution of a stochastic integral
equation with degeneration\/}, 
Izvestiya Akademii Nauk SSSR, seriya matematicheskaya,
Vol. 36 (1972), No. 1,  248--261 in Russian; English translation:
Math. USSR
Izvestija, Vol. 6 (1972),  No. 1, 249--262.

\bibitem{Kr97} N.V. Krylov,  {\em On the rate of convergence
of finite--difference approximations for Bellman's equations\/},
 Algebra i Analiz,  
 Vol. 9 (1997), No. 3, 245--256 in Russian;
English translation: St. Petersburg Math. J.,
Vol. 9 (1998), No. 3, 639--650.

\bibitem{Kr99} N.V. Krylov, {\em Approximating
value functions for  controlled degenerate diffusion
processes by using piece-wise constant policies},
Electronic Journal of Probability, Vol.4 (1999), paper no. 2, 1--19,\\
http://www.math.washington.edu/\~\,ejpecp/EjpVol4/paper2.abs.html

\bibitem{Kr07} N.V. Krylov, {\em A priori estimates
of smoothness of solutions to  difference  Bellman  equations 
  with linear and quasi-linear operators\/}, 
Math. Comp., Vol. 76 (2007), 669-698. 

\bibitem{Kr08} N.V. Krylov, {\em On
factorizations of smooth nonnegative matrix-values
functions and on
smooth functions with values in polyhedra\/},
  Appl. Math. Optim., Vol. 58 (2008), No. 3, 373--392.

\bibitem{Kr11} N.V. Krylov,  {\em On a representation of fully nonlinear elliptic
operators in terms of pure
second order derivatives and its applications\/}, 
Problems in Mathematical Analysis, Vol.
 59, July 2011, p. 3--24 in Russian,
English translation in Journal of
Mathematical Sciences (New York),
DOI: 10.1007/s10958-011-0445-0.

\bibitem{Kr12.1} N.V. Krylov, {\em Interior estimates for 
the first-order differences
for finite-difference
  approximations for elliptic Bellman's
equations\/}, to appear in Applied Math. Optimiz.


\bibitem{Kr12.2} N.V. Krylov, {\em Interior estimates for second differences 
  of solutions of finite-difference elliptic Bellman's 
equations\/}, to appear in Math. Comp.

\bibitem{Ku} H.J. Kushner, ``Probability methods for 
approximations
in stochastic control and for elliptic equations'', 
Mathematics in Science
and Engineering, Vol. 129. Academic Press [Harcourt 
Brace Jovanovich,
Publishers], New York--London, 1977.

\bibitem{KD} Kushner  H.J. and
 Dupuis  P.G., ``Numerical methods
for stochastic control problems
in continuous time'',
Springer Verlag, 1992.

\bibitem{KT90} H.-J,  Kuo and N.S. Trudinger,
{\em Linear elliptic difference inequalities with random coefficients\/},
Math. Comp., Vol. 55 (1990), No. 191, 37--53.

\bibitem{KT92} H.-J,  Kuo and N.S. Trudinger, {\em
Discrete methods for fully nonlinear elliptic equations\/},
  SIAM Journal on Numerical Analysis, Vol. 29 (1992), No. 1,  123--135.

 
\bibitem{KT96} H.-J,  Kuo and N.S. Trudinger, {\em Positive 
difference operators on general meshes}, Duke Math. J., Vol. 83
(1996), No. 2, 415--433.

\bibitem{Pr} G. Pragarauskas, 
{\em Approximation of controlled solutions of
It\^o equations by controlled 
Markov chains\/}, Lit. Mat. Sbornik, Vol. 
23 (1983), No. 1, 175--188 in Russian;
English translation: Lithuanian
 Math. J., Vol. 23 (1983), No. 1, 98--108.

\bibitem{So} P.E. Souganidis, {\em Approximation schemes for viscosity 
solutions of Hamilton-Jacobi equations\/},
J. Differential Equations, Vol. 59 (1985), No. 1, 1--43.


\bibitem{Sa} M.V. Safonov, {\em
Classical solution of second-order nonlinear elliptic 
equations\/},  Izv. Akad.
Nauk SSSR Ser. Mat., Vol. 52  (1988),  No. 6, 1272--1287  in Russian;
  English translation:  Math.
USSR-Izv., Vol.  33  (1989),  No. 3, 597--612.

\end{thebibliography}
\end{document}